\documentclass[a4paper,12pt]{article}
\usepackage{amsfonts}
\usepackage{amsmath}
\usepackage{amssymb}
\usepackage{amsthm}
\usepackage{diagrams}
\usepackage{hyperref}
\usepackage[pdftex]{graphicx}

\newtheorem{theorem}[equation]{Theorem}
\newtheorem{corollary}[equation]{Corollary}
\newtheorem{lemma}[equation]{Lemma}
\newtheorem{proposition}[equation]{Proposition}

\numberwithin{equation}{section}

\newcommand{\Z}{\mathbb{Z}}

\newcommand{\diag}{\operatorname{diag}}
\newcommand{\bq}{/\!\!/}

\begin{document}

\title{Three new almost positively curved manifolds}

\author{Jason DeVito}

\date{}

\maketitle

\begin{abstract}

A Riemannian manifold is called almost positively curved if the set of points for which all $2$-planes have positive sectional curvature is open and dense.  We  find three new examples of almost positively curved manifolds:  $Sp(3)/Sp(1)^2$, and two circle quotients of $Sp(3)/Sp(1)^2$.  We also show the quasi-positively curved metric of Tapp \cite{Ta1} on $Sp(n+1)/Sp(n-1) Sp(1)$ is not almost positively curved if $n\geq 3$.
\end{abstract}

\section{Introduction}

There are few known examples of simply connected compact manifolds admitting metrics of positive sectional curvature.  In fact, other than the rank one symmetric spaces, there are only two infinite families known (one in dimension 7 \cite{AW,Es1}, the other in 13 \cite{Baz1,Ber}) and seven other examples: two in dimension $6$, two in dimension $7$, and single examples in dimension $12$ and $24$.\cite{Es1,Wa,Ber,GVZ,De}.

Despite the paucity of examples, there are few obstructions distinguishing compact manifolds admitting metrics of non-negative curvature from those admitting metrics with positive curvature.  Further, all known obstructions vanish on compact simply connected manifolds, that is, there is no known example of a compact simply connected non-negatively curved Riemannian manifold which does not admit a positively curved metric.

As a means to understand the difference between manifolds admitting non-negative curvature from those admitting positive curvature, much work has gone into understanding two classes of Riemannian manifolds which lie in between- the quasi-positively curved manifolds and the almost positively curved manifolds.  
Recall that a Riemannian manifold is said to be \textit{quasi-positively curved} if it is everywhere non-negatively curved and has a point $p$ (and therefore an open neighborhood) for which the sectional curvatures of all 2-planes at $p$ are positive.  A manifold is said to be \textit{almost positively curved} if it has an open dense subset of points for which all sectional curvatures are positive.

Examples with these weaker curvature notions are much more abundant than in the case of strictly positive curvature \cite{Wi,Ta1,KT,DDRW,DM,DN,Ke1,Ke2,EK,GrMe1,PW2,W,Di1}.  The main constructions are due to Wilking \cite{Wi} and Tapp \cite{Ta1}.  In \cite{KT}, Kerr and Tapp classify so-called positive triples, spaces for which the metrics constructed in \cite{Ta1} are quasi-positively curved with a point of positive curvature arbitrarily close to the identity coset.  In particular, they show that $M_n:=Sp(n+1)/Sp(n-1)Sp(1)$ admits a Riemannian metric of quasi-positive curvature.  Here, $Sp(1)\times Sp(n-1)$ is embedded into $Sp(n+1)$ via the block embedding.  The metric Kerr and Tapp use admits two different free isometric $Sp(1)$ actions.  Quotienting by one gives the homogeneous space $Sp(n+1)/Sp(n-1)Sp(1)^2$, and quotienting by the other gives the biquotient $\Delta Sp(1)\backslash Sp(n+1)/Sp(n)Sp(1)$, where $\Delta Sp(1)$ indicates the block embedding $Sp(1)\rightarrow Sp(n+1)$ with $q\mapsto \diag(q,q,...,q)$.  Wilking \cite{Wi} has shown both of these $Sp(1)$ quotients of $M_n$ admit metrics of almost positive curvature.  When $n=2$, the homogeneous space $Sp(3)/Sp(1)^3$ admits a homogeneous metric with strictly positive curvature \cite{Wa}.

If we restrict the above $Sp(1)$ actions to $S^1\subseteq Sp(1)$, we obtain manifolds $Q_n:=Sp(n+1)/Sp(n-1)Sp(1)S^1$ and $R_n := \Delta S^1 \backslash Sp(n+1)/Sp(n-1)Sp(1)$.  It immediately follows from O'Neill's formulas \cite{On1} that each $Q_n$ and $R_n$ admits a metric of quasi-positive curvature.  Thus, one is naturally led to wonder which of the $M_n, Q_n,$ and $R_n$ admit metrics of almost positive curvature.

\begin{theorem}\label{main}  The homogeneous spaces $Sp(3)/Sp(1)^2$, $Sp(3)/Sp(1)^2S^1$ and the biquotient $\Delta S^1 \backslash Sp(3)/Sp(1)^2$ each admit Riemannian metrics of almost positive curvature.

\end{theorem}

The cohomology rings of $Q_2$ and $R_2$ are isomorphic, but their Pontryagin classes distinguish them up to homotopy, see Propositions \ref{cohomring} and \ref{homotopytype}.  The space $M_2$ is a parallelizable $S^4$ bundle over $S^{11}$, with cohomology ring isomorphic to that of $S^4\times S^{11}$.  Nonetheless, we show it is not even homotopy equivalent to $S^4\times S^{11}$, being distinguished by their $10$th homotopy groups (Proposition \ref{nothe}).

In contrast with Theorem \ref{main}, we show that Tapp's metrics are not always almost positively curved.

\begin{theorem}\label{main2} Tapp's metrics on $Sp(n+1)/Sp(n-1)Sp(1)$ for $n\geq 3$ are not almost positively curved.

\end{theorem}

Instead of working with Tapp's construction, we use Wilking's construction \cite{Wi} which gives, up to scale, an isometric metric \cite{Ta1}.  That is, we express each of the manifolds as a biquotient of the form $\Delta G\backslash G\times G/U$.  Beginning with a bi-invariant metric on $G\times G$, we Cheeger deform each metric to find a new non-negatively curved left invariant metric on $G\times G$ which has fewer zero-curvature planes, but for which $U$ still acts isometrically.  The metric on $\Delta G\backslash G\times G/U$ is the induced metric coming from the canonical submersion $G\times G\rightarrow \Delta G\backslash G\times G/U$.

Our main tool for proving Theorems \ref{main} and \ref{main2} is identifying an explicit two-dimensional disk in $M_n$ which intersects every orbit of the isometry group of $M_n$, see Proposition \ref{isom} and Lemma \ref{lem:so3form}.  This reduces the computations to a two-dimensional set of points, a considerable simplification.

The outline of the paper is as follows.  In Section 2, after covering the necessary background information, we define the metrics on $M_n$ and prove that the action of the isometry group is of cohomogeneity two, that is, that the quotient $M_n/Iso(M_n)$ is two-dimensional.  In Section \ref{lowdimalmpos}, we prove Theorem \ref{main}.  In Section \ref{hidim}, we determine an open set of points of $M_n$ where each point has an infinite number of distinct zero-curvature planes, proving Theorem \ref{main2}.  In Section \ref{top}, we compute the topology of $M_2, Q_2$, and $R_2$, showing that these examples are distinct, up to homotopy, from any previously known example with almost positive curvature.

\

Theorem \ref{main} was originally proven in the author's thesis.  He is greatly indebted to Wolfgang Ziller for helpful comments.

\section{Biquotients and their geometry}\label{geo}

In this section, we first recall the relevant information about the geometry of biquotients, and then apply it to our specific examples.

\subsection{Background}

Given a compact Lie group $G$ and a subgroup $U\subseteq G\times G$, there is a natural action of $U$ on $G$ given by $(u_1,u_2)*g = u_1 g u_2^{-1}$.  The action is effectively free iff whenever $u_1$ and $u_2$ are conjugate, we have $u_1=u_2\in Z(G)$.  When the action is effectively free, the orbit space $G\bq U$ naturally has the structure of a manifold such that $\pi:G\rightarrow G\bq U$ is a submersion.  The orbit space is called a biquotient.  When $U\subseteq \{e\}\times G\subseteq G\times G$, the action is automatically free and the quotient is the homogeneous space $G/U$.

Suppose $U$ acts on $G$ effectively freely.  If $G$ is equipped with any $U$-invariant metric, then $G\bq U$ inherits a unique Riemannian metric for which the projection is a Riemannian submersion.  We will sometimes refer to this Riemannian metric on $G\bq U$ as the metric induced from the Riemannian metric on $G$.   By equipping $G$ with a bi-invariant metric, we see, via O'Neill's formulas \cite{On1}, that every biquotient admits a metric of non-negative sectional curvature.

To find metrics on $G\bq U$ with fewer zero-curvature planes, two main techniques are used:  Cheeger deformations and Wilking's doubling trick.  We now focus on Cheeger deformations.

The idea of a Cheeger deformation is to shrink a left $G$-invariant Riemannian metric on Lie group in the direction of a subgroup.  If the original Riemannian metric is non-negatively curved, then the deformed metric is also non-negatively curved.  In addition, the deformed metric tends to have fewer zero-curvature planes than the original metric, but it also tends to have a smaller isometry group.  The precise versions of these statements can be found in \cite{Es1}; we describe them  below in order to set up notation.

Let $K$ be a closed subgroup of $G$ and equip $G$ with left $G$-invariant, right $K$-invariant Riemannian metric $\langle \cdot, \cdot \rangle_0$ of non-negative sectional curvature.  For example, $\langle \cdot, \cdot \rangle_0$ could be a bi-invariant metric on $G$.  For each $t > 0$, we equip $G\times K$ with the product Riemannian metric $\langle \cdot, \cdot \rangle_0 + t\langle \cdot, \cdot \rangle_0|_K$.  Then, for each $t$, $K$ acts freely and isometrically on $G\times K$ via $k\ast (g_1,k_1) = (g_1k^{-1}, kk_1)$.  The quotient $G\times_K K$, which is diffeomorphic to $G$ via the diffeomorphism $[(g_1,k_1)]\mapsto g_1k_1$, thus inherits a unique Riemannian metric for which the submersion $G\times K\rightarrow G\times_K K\cong G$ is a Riemannian submersion.  We denote this induced Riemannian metric on $G$ by $\langle \cdot, \cdot \rangle_1$, and we call it the Cheeger deformation of $\langle \cdot, \cdot \rangle_0$ in the direction of $K$.  We note that the metric on $G\times K$ is a product of non-negatively curved metrics, so is non-negatively curved.  It follows from O'Neill's formulas for a submersion \cite{On1} that $\langle \cdot, \cdot \rangle_1$ is also non-negatively curved.

The isometry group of $\langle \cdot,\cdot \rangle_1$ is large, though, in general, it is smaller than that of $\langle \cdot , \cdot \rangle_0$.  Specifically, the action of $G\times K$ on itself given by $(g,k)\ast (g_1, k_1) = (gg_1, k_1 k^{-1})$ is isometric and commutes with the free $K$ action above, so descends to an isometric action on $G\times_K K \cong G$.  In particular, $\langle \cdot, \cdot \rangle_1$ is a left $G$-invariant, right $K$-invariant non-negatively curved metric.  As the left $G$-action is transitive, we can reduce all curvature computations to a single point, the identity  $e\in G$.

Writing $\mathfrak{g} = T_e G$ for the Lie algebra of $G$ and $\mathfrak{k}$ for the Lie algebra of $K$, we get a decomposition $\mathfrak{g} = \mathfrak{k}\oplus \mathfrak{p}$, orthogonal with respect to $\langle \cdot, \cdot \rangle_0$.  Then, one can show (see e.g. \cite{EK}) that for $X,Y\in \mathfrak{g}$, $\langle X,Y\rangle_1 = \langle\phi(X), Y\rangle_0$, where $\phi:\mathfrak{g}\rightarrow\mathfrak{g}$ is defined by $\phi(X) = X_\mathfrak{p} + \frac{t}{t+1}X_\mathfrak{k}$.

Writing $\sec_1$ for sectional curvature with respect to the metric $\langle \cdot, \cdot \rangle_1$, Eschenburg \cite{es2} proved the following proposition.

\begin{proposition}\label{1def}  If $\langle \cdot, \cdot \rangle_0$ is a bi-invariant metric and $(G,K)$ is a symmetric pair, then $\sec_1(\phi^{-1} X, \phi^{-1}Y) = 0$ iff $[X,Y] = [X_\mathfrak{k}, Y_\mathfrak{k}] = [X_{\mathfrak{p}}, Y_{\mathfrak{p}}] = 0$.

\end{proposition}

We recall that for a bi-invariant metric $\langle \cdot, \cdot\rangle_0$, that $\sec_0(X,Y) = 0$ iff $[X,Y]=0$.  Thus, Proposition \ref{1def} indicates that,  in general, $\langle \cdot , \cdot \rangle_1$ has fewer zero-curvature planes than $\langle \cdot,\cdot\rangle_0$.

\

Wilking's doubling trick is based up an observation of Eschenburg, that the map $\psi:G\times G\rightarrow G$ given by $\psi(g_1,g_2) = g_1^{-1} g_2$ induces a diffeomorphism $\overline{\psi}:\Delta G \backslash G\times G/U\rightarrow G\bq U$, where $\Delta G\times U$ acts on $G\times G$ via $(g,(u_1,u_2))\ast(g_1,g_2) = (gg_1 u_1^{-1}, gg_2 u_2^{-1})$.  Wilking \cite{Wi} noticed that one can choose a Cheeger deformed Riemannian metric on each factor, giving a larger class of natural metrics of non-negative curvature.  Tapp \cite{Ta1} has shown that his metrics are, up to scale, isometric to those of Wilking's where one uses the same metric on both factors.

Suppose $\langle \cdot, \cdot \rangle_1$ is obtained from a bi-invariant metric $\langle \cdot, \cdot\rangle_0$ on $G$ by Cheeger deforming in the direction of $K\subseteq G$.  If $U\subseteq K\times K\subseteq G$, then $U$ acts isometrically on $G\times G$ equipped with the product metric $\langle \cdot, \cdot \rangle_1 + \langle \cdot, \cdot \rangle_1$; let $\langle \cdot, \cdot \rangle_2$ denote the induced metric on $G\bq U$, that is $\langle\cdot,\cdot\rangle_2$ is the unique Riemannian metric on $G\bq U$ for which the submersion $G\times G\rightarrow \Delta G\backslash (G\times G)/U\cong G\bq U$ is a Riemannian submersion.  We wish to understand when a 2-plane $\sigma$ in $G\bq U$ has zero sectional curvature with respect to $\langle \cdot, \cdot \rangle_2$.  By O'Neill's formula, if $\sigma$ has zero sectional curvature with respect to $\langle \cdot,\cdot\rangle_2$, then the horizontal lift of $\sigma$ must have zero-curvature.  Thus, we must determine the horizontal distribution on $G\times G$ with respect to $\langle \cdot, \cdot \rangle_1 + \langle \cdot, \cdot \rangle_1$.

It is clear that every orbit of the $\Delta G\times U$ action passes through a point of the form $(g_1,e)$, where $e\in G$ is the identity element, so we may focus on determining the horizontal space at points of this form.  Further, since $\langle \cdot, \cdot \rangle_1 + \langle\cdot,\cdot\rangle_1$ is left $G\times G$-invariant, the curvature of a $2$-plane at $T_{(g_1,e)} G\times G$ is the same as the curvature of its left translations.  Thus, instead of directly computing the vertical space $\mathcal{V}_{g_1}$ and horizontal space $\mathcal{H}_{g_1}$ at the point $(g_1,e)$, we use the left translations of these subspaces to $(e,e)\in G\times G$.

As is shown in, e.g. \cite{Ke1}, the vertical subspace $\mathcal{V}_{g_1}$ at $(g_1,e)\in G\times G$, translated to $(e,e)$ using left translation, is $$ (L_{g_{1}^{-1}})_\ast\mathcal{V}_{g_1} = \{((Ad_{g_{1}^{-1}}X)-U_1, X-U_2)|\; X\in\mathfrak{g} \text{ and } (U_1,U_2)\in\mathfrak{u}\}$$ where $\mathfrak{u}\subseteq \mathfrak{g}\oplus \mathfrak{g}$ is the Lie algebra of $U$.

\begin{proposition}\label{horiz} With respect to $\langle \cdot, \cdot \rangle_1$, the horizontal space $(L_{g_{1}^{-1}})_\ast \mathcal{H}_{g_1}$ at the point $(g_1,e)$, left translated to $(e,e)$, is \begin{align*} \{(\phi^{-1}(-Ad_{g_{1}^{-1}} X),\phi^{-1}(X)):\\ X\in \mathfrak{g} \text{ and }&  \langle X,Ad_{g_1} U_1 - U_2\rangle_0 = 0 \text{ for all } (U_1,U_2)\in\mathfrak{u}\}.\end{align*}\end{proposition}

\begin{proof}
Consider the linear map $\rho_{g_1}:\mathfrak{u}\rightarrow \mathfrak{g}$ with $\rho_{g_1}(U_1,U_2) = Ad_{g_1}U_1 - U_2$.  If $(U_1,U_2)\in \ker \rho_{g_1}$, then $(U_1,U_2) = (U_1, Ad_{g_1}U_1) \in \mathfrak{u}$, so $(u_1, g_1 u_1 g_1^{-1}) := \exp(U_1, Ad_{g_1} U_1)\in U$.  Then, recalling the biquotient action of $U$ on $G$, we see $(u_1 ,g_1 u_1 g_1^{-1})\ast g_1^{-1} = u_1 g_1^{-1} (g_1 u_1 g_1^{-1})^{-1} = g_1^{-1}$, that is, $(u_1, g_1 u_1 g_1^{-1})\in U$ fixes $g_1^{-1}\in G$.  Since the $U$ biquotient action on $G$ is free, we conclude that $u_1 = \exp(U_1) = e$.  Repeating this argument for every real multiple of $(U_1,U_2)$, we deduce that $\exp(sU_1) = e$ for all $s\in \mathbb{R}$.  Thus, $U_1 = 0$ and $U_2 = Ad_{g_1}U_1 = 0$ as well.  That is, $\ker \rho_{g_1} = \{0\}$.

It follows that $(\rho_{g_1}(\mathfrak{u}))^\bot = \{X \in \mathfrak{g}: \langle X, \rho_{g_1}(\mathfrak{u})\rangle_0 = 0\}$ has dimension $\dim G - \dim U$, the same $\dim (L_{g_{1}^{-1}})_\ast \mathcal{H}_{g_1}$.  So, to establish the proposition, it is sufficient to show each $(\phi^{-1}(-Ad_{g_1^{-1}} X), \phi^{-1}(X))$ with $X\in (\rho_{g_1}(\mathfrak{u}))^\bot$ is an element of $(L_{g_{1}^{-1}})_\ast \mathcal{H}_{g_1}$.

So, consider the element $(\phi^{-1}(-Ad_{g_1^{-1}} X), \phi^{-1}(X))$ with $X\in (\rho_{g_1}(\mathfrak{u}))^\bot$.  Because $Ad_{g_1}$ is an isometry of $\langle \cdot, \cdot \rangle_0$, we compute  \begin{align*} \langle  \phi^{-1}(-Ad_{g_1^{-1}} X), Ad_{g_1^{-1}} Y\rangle_1 + \langle \phi^{-1}(X), Y\rangle_1 &= \langle -Ad_{g_1^{-1}} X, Ad_{g_1^{-1}} Y\rangle_0 + \langle X,Y\rangle_0\\ &= -\langle X,Y\rangle_0 + \langle X,Y\rangle_0\\ &= 0  \end{align*}

Similarly, we see \begin{align*} \langle \phi^{-1}(-Ad_{g_1^{-1}} X), U_1\rangle_1 + \langle \phi^{-1} X, U_2\rangle_1 &= -\langle Ad_{g_1^{-1}} X, U_1\rangle_0 + \langle  X, U_2\rangle_0\\ &=-\langle X, Ad_{g_1} U_1\rangle_0 + \langle X, U_2\rangle_0\\ &= 0\end{align*} since $X\in (\rho_{g_1}(\mathfrak{u}))^\bot$.

Thus, we conclude $(\phi^{-1}(-Ad_{g_1^{-1}} X), \phi^{-1}(X))\in (L_{g_1^{-1}})_\ast \mathcal{H}_{g_1}$.  The proposition now follows.   
\end{proof}

Since $\langle \cdot, \cdot \rangle_1 + \langle \cdot, \cdot \rangle_1$ is a product of non-negatively curved metrics, we obtain the following corollary.

\begin{corollary}\label{doublezerocurv} Suppose $\operatorname{span}\{(\phi^{-1}(-Ad_{{g_1}^{-1}} X), \phi^{-1}X),(\phi^{-1}(-Ad_{{g_1}^{-1}} Y), \phi^{-1}Y)\}$ is horizontal and has zero-curvature.  Then the $2$-planes $$\operatorname{span}\{\phi^{-1}(Ad_{g_1}^{-1} X), \phi^{-1}(Ad_{g_1}^{-1}Y)\} \text{ and } \operatorname{span}\{ \phi^{-1}(X), \phi^{-1}(Y)\} $$ both have zero-curvature with respect to $\langle \cdot, \cdot \rangle_1$.

\end{corollary}
 
We now assume $(G,K)$ is a symmetric pair, so Proposition \ref{1def} applies to both planes given in Corollary \ref{doublezerocurv}.  In particular both $$[Ad_{g_1^{-1}} X, Ad_{g_1^{-1}} Y] = [(Ad_{g_1^{-1}} X)_{\mathfrak{k}}, (Ad_{g_1^{-1}} Y)_{\mathfrak{k}}] =  [(Ad_{g_1^{-1}} X)_{\mathfrak{p}}, (Ad_{g_1^{-1}} Y)_{\mathfrak{p}}] = 0$$ and $$[X,Y] = [X_{\mathfrak{k}}, Y_{\mathfrak{k}}] = [X_{\mathfrak{p}}, Y_{\mathfrak{p}}] = 0.$$

However, some of these conditions are redundant.  Indeed, since $(G,K)$ is symmetric, then with respect to an $\langle \cdot, \cdot \rangle_0$-orthogonal decomposition $\mathfrak{g} = \mathfrak{k}\oplus \mathfrak{p}$, we have $[\mathfrak{p}, \mathfrak{p}]\subseteq \mathfrak{k}$.  Then, under the assumption $[X,Y] = 0$, it follows easily that $[X_\mathfrak{k}, Y_\mathfrak{k}] = 0$ iff $[X_\mathfrak{p}, Y_\mathfrak{p}] = 0$.  Further, since $Ad_{g_1}$ is a Lie algebra isomorphism, $[X,Y] = 0$ iff $[Ad_{g_1} X, Ad_{g_1}Y] = 0$.

It follows that the $2$-planes $$\operatorname{span}\{\phi^{-1}(Ad_{g_1}^{-1} X), \phi^{-1}(Ad_{g_1}^{-1}Y)\} \text{ and } \operatorname{span}\{ \phi^{-1}(X), \phi^{-1}(Y)\} $$ both have zero-curvature with respect to $\langle\cdot, \cdot\rangle_1$ iff $$[X,Y] = [X_\mathfrak{p}, Y_\mathfrak{p}] = [(Ad_{g_1}^{-1}X)_\mathfrak{p},(Ad_{g_1}^{-1}Y)_\mathfrak{p}] = 0.$$  We note that these conditions on $X$ and $Y$ really only depend on $\operatorname{span}\{X,Y\}$.

With nicer assumptions on $U$, $K$, and $G$, even more is true.

\begin{proposition}\label{eqns}

Suppose $(G,K)$ is a symmetric pair with $U\subseteq K$ and suppose that a bi-invariant metric $\langle\cdot, \cdot \rangle_0$ on $G$ induces a positively curved metric on $G/K$ and that its restriction $\langle \cdot, \cdot \rangle|_K$ induces a positively curved metric on $K/U$.  Suppose $\langle \cdot, \cdot \rangle_1$ is obtained by Cheeger deforming $\langle \cdot, \cdot \rangle_0$ in the direction of $K$.  Then, with respect to the Riemannian submersion metric $\langle \cdot, \cdot\rangle_2$ on $G/U\cong \Delta G\backslash (G\times G, \langle \cdot,\cdot\rangle_1 + \langle \cdot, \cdot \rangle_1)/U$, there is a zero-curvature plane at the point $[(g_1,e)]\in \Delta G\backslash G\times G/ U$ iff there are non-zero vectors $X = X_{\mathfrak{k}},Y = Y_{\mathfrak{p}} \in \mathfrak{g}$ satisfying each of the following three conditions.

\begin{equation}
\langle X,\mathfrak{u}\rangle_0  = 0 \label{con1} \tag{Condition A}
\end{equation}

\begin{equation}
[X,Y] = 0 \label{con2} \tag{Condition B}
\end{equation}

\begin{equation}
(Ad_{g_1^{-1}} X)_{\mathfrak{p}} \text{ and } (Ad_{g_1^{-1}} Y)_{\mathfrak{p}} \text{ are dependent over }\mathbb{R}. \label{con3} \tag{Condition C}
\end{equation}

\end{proposition}

\begin{proof}
First, assume there is a zero-curvature plane at $[(g_1,e)]\in G/U$.  From O'Neill's formula \cite{On1}, there must be a horizontal zero-curvature plane at $(g_1,e)\in G\times G$.  By Proposition \ref{horiz}, Corollary \ref{doublezerocurv}, and the following discussion, there are vectors $$(\phi^{-1}(-Ad_{g_{1}^{-1}} X),\phi^{-1}(X)) \text{ and } (\phi^{-1}(-Ad_{g_{1}^{-1}} Y),\phi^{-1}(Y))\in (L_{g_1^{-1}})_\ast \mathcal{H}_{g_1}$$ (so both $X$ and $Y$ are $\langle \cdot, \cdot\rangle_0-$orthogonal to $\mathfrak{u}$) which satisfy $$[X,Y] = [X_\mathfrak{p}, Y_\mathfrak{p}] = [(Ad_{g_1^{-1} }X)_\mathfrak{p},(Ad_{g_1^{-1} }Y)_\mathfrak{p}] = 0.$$  In particular, \ref{con1} and \ref{con2} follow.

Recall that the curvature of a normal homogeneous space $G/K$ is given by $\sec(X,Y) = \frac{1}{4}\| [X,Y]_\mathfrak{p}\|^2 + \| [X,Y]_\mathfrak{k}\|^2$ with $\mathfrak{g} = \mathfrak{k}\oplus \mathfrak{p}$ (see, e.g. \cite[Corollary 3.33]{CE}).  In particular, $\sec(X,Y) = 0$ iff $[X,Y]=0$.

Applying this to $X_\mathfrak{p}$ and $Y_\mathfrak{p}$, interpreted as elements of $T_{eK} G/K$, we see that $[X_\mathfrak{p}, Y_\mathfrak{p}] = 0$ iff $X_\mathfrak{p}$ and $Y_\mathfrak{p}$ are linearly dependent over $\mathbb{R}$.  This same argument also applies to $(Ad_{g_1^{-1}}X)_\mathfrak{p}$ and $(Ad_{g_1^{-1}} Y)_\mathfrak{p}$, giving \ref{con3}.  By subtracting an appropriate multiple of $Y$ from $X$, we find a vector $X' = X'_\mathfrak{k}$ for which $\operatorname{span}\{X,Y\} = \operatorname{span}\{X',Y\}$.

Similarly, since both $X'$ and $Y$ are $\langle \cdot, \cdot\rangle_0-$orthogonal to $\mathfrak{u}$, we may interpret $X' = X'_\mathfrak{k}, Y_\mathfrak{k}$ as elements of $T_{eU} K/U$.  In particular, $[X', Y_\mathfrak{k}] = 0$ iff $X'$ and $Y_\mathfrak{k}$ are linearly dependent over $\mathbb{R}$.  Then, by subtracting an appropriate multiple of $X'$ from $Y$, we find a new vector $Y' = Y'_\mathfrak{p}$ for which $\operatorname{span}\{X,Y\} = \operatorname{span}\{X',Y'\}$.  Then $X'\in\mathfrak{k}$, $Y'\in \mathfrak{p}$, and $X'$ and $Y'$ satisfy \ref{con1}, \ref{con2}, and \ref{con3}.

\

Conversely, assume there is an $X = X_\mathfrak{k}$ and $Y = Y_\mathfrak{p}$ satisfying all three conditions.  Note that $Y$ is automatically $\langle \cdot, \cdot\rangle_0-$orthogonal to $\mathfrak{u}$ because $\mathfrak{u}\subseteq \mathfrak{k}$ and $\mathfrak{p}$ is $\langle \cdot, \cdot\rangle_0-$orthogonal to $\mathfrak{k}$.

It follows that the vectors $$(\phi^{-1}(-Ad_{g_{1}^{-1}} X),\phi^{-1}(X)) \text{ and } (\phi^{-1}(-Ad_{g_{1}^{-1}} Y),\phi^{-1}(Y))$$ are elements of $\mathcal{H}_{g_1}$ and from Proposition \ref{1def}, that they span a horizontal zero-curvature plane in $G\times G$.  Finally, Tapp \cite{Ta2} has shown that in this setup, a horizontal zero curvature plane projects to a zero-curvature plane in $G/U$.  

\end{proof}

All of the hypothesis of Proposition \ref{eqns} will apply to the $M_n$ examples.

\subsection{Applications to our examples}\label{metrics}

In this section, we apply the discussion in the previous section to define the metrics on $M_n = Sp(n+1)/Sp(n-1)Sp(1)$.  We also show that the isometry group acts by cohomogeneity two, finding a nice section for the action.

We let $G = Sp(n+1)$ denote the group of $(n+1)\times (n+1)$ unitary matrices over $\mathbb{H}$.  We let $K = Sp(n)\times Sp(1)$ block embedded into $G$ and we set $U\cong Sp(n-1)\times Sp(1)$, embedded in $G$ via $(A,q)\in U\mapsto \diag(A, 1, q)$.  We will also consider the subgroup $N\subseteq K$ with $N = \{\diag(A,q_1,q_2):(A,q_2)\in U, q_1\in Sp(1)\}$.  Note that $N$ normalizes $U$, so $N/U\cong Sp(1)$ acts on $G/U$, and this action is isometric as $N\subseteq K$.  However, we stress that $N$ does not normalize $Sp(n-1)S^1 Sp(1)$, and $K$ does not normalize $\Delta S^1$, so the following arguments do not apply to the two circle quotients $Q_n$ and $R_n$ of $M_n$.

The inner product on $\mathfrak{g}$ given by $\langle X,Y\rangle = -\operatorname{Re}Tr(XY)$ is $Ad_G$ invariant, so extends to a unique bi-invariant Riemannian metric $\langle \cdot,\cdot\rangle_0$ on $G$.  We Cheeger deform $\langle \cdot, \cdot\rangle_0$ in the direction of $K$ and call the resulting metric $\langle\cdot,\cdot\rangle_1$.  Equipping $G\times G$ with the product metric $\langle \cdot, \cdot \rangle_1 + \langle \cdot, \cdot \rangle_1 $, the natural action by $\Delta G \times K\times N$ given by $(g,k,n)\ast (g_1, g_2) = (g \, g_1\, k^{-1}, g\, g_2 n^{-1})$ is isometric and the restriction of the action to $\Delta G \times \{1\} \times U$ is free.  We give the quotient $\Delta G \backslash G\times G/U$, which is canonically diffeomorphic to $G/U$, the submersion metric $\langle \cdot, \cdot\rangle_2$, as in the previous section.

We now show the isometry group acts with a two dimensional quotient space.  To do so, we use the following notation.  For any $g\in G = Sp(n+1)$, we let $v(g)$ denote the last column of $g$, interpreted as an element of $S^{4n+3}\subseteq \mathbb{H}^{n+1}$.  We also use $v_0(g)$ to denote the first $n-1$ entries of $v(g)$, $v_n(g)$ to denote the second to last entry of $v(g)$, and $v_{n+1}(g)$ to denote the last entry of $v(g)$.

\begin{proposition}\label{isom}  

Consider the $\Delta G\times K\times N$ action on $G\times G$ given by $(g,k,n)\ast(g_1,g_2) = (gg_1 k^{-1}, gg_2 n^{-1})$.  Then $(g_1,g_2)$ and $(h_1,h_2)$ are in the same orbit iff $|v_i(g_2^{-1} g_1)| = |v_i(h_2^{-1} h_1)|$ for each of $i = 0,n,n+1$.

\end{proposition}

\begin{proof}

Let $(g_1, g_2)\in G\times G$.  For any $h\in N$, we set $g = h g_2^{-1}$, so $(g,k,h)\ast (g_1,g_2) = (h \, g_2^{-1} \, g_1 \, k^{-1} , e)$.  Thus, we need only show that $g_1, g_2\in G$ are equivalent under the $N\times K$ on $G$ given by $(h,k) \ast g = hg k^{-1}$ iff $|v_i(g_1)| = |v_i(g_2)|$ for $i = 0, n, n+1$.

Now, let $g\in G$.  Then, for $h = \diag(A,q_1, q_2)\in N$, it is easy to verify that $v_0(hg) = Av_0(g)$, $v_n(hg) = q_1 v_n(g)$, and $v_{n+1}(hg) = q_2 v_{n+1}(g)$.  Since left multiplication by elements of $Sp(n-1)$ and $Sp(1)$ preserves lengths, it now follows that $|v_i(hg)| = |v_i(g)|$ for $i = 0,n,n+1$.

Likewise, for an element $k = \diag(B,q)\in K = Sp(n)\times Sp(1)$, we have $v_i(g k^{-1}) = v_i(g) q^{-1}$, so $|v_i(gk^{-1})| = |v_i(g)|$ for $i=0,n,n+1$.  This establishes the fact that the $N\times K$ action preserves each $|v_i(g)|$.  We now show that these are the only invariants.

To that end, first note that for $Sp(n)\times \{1\}\subseteq K$,  $G/Sp(n) \cong S^{4n+3}$, with the diffeomorphism induced by mapping $g\in G$ to $v(g)$.  Thus, the orbit through $g$ is determined by $v(g)$.  In fact, since left multiplication on $Sp(n-1)$ (resp. $Sp(1)$) is transitive on the unit sphere in $\mathbb{H}^{n-1}$ (resp. $\mathbb{H}$), for each $g\in G$, there is an $h = \diag(A,q_1,q_2)\in N$ for which $v(hg)$ has entries which are all zero, except for the last three which are the non-negative real numbers $|v_i(hg)|$ for $i=0,n,n+1$.  This shows that the orbit through $g$ is completely determined by $|v_i(g)|$ for $i = 0,n,n+1$.   
\end{proof}

Given any column vector $w=(w_1,w_2,w_3)^t\in \mathbb{R}^3$ with unit length, there is a matrix $A \in SO(3)$ for which $A_{12} = 0$ and for which the last column of $A$ is $w$.  To see this, note that $w^\bot \cap \{(0,x_2,x_3)^t\}\subseteq \mathbb{R}^3$ has dimension at least one, so we can pick a non-zero element of the intersection to use as the second column of $A$.  In fact, we can pick this non-zero vector to have a non-negative second entry.  The first column must then be the cross product of the second and third.  

We let $$\mathcal{F} =\{A\in SO(3): A_{12} = 0, A_{ij}\geq 0\text{ for } (i,j)\in\{(2,2),(1,3),(2,3),(3,3)\}\}.$$  We will identify $\mathcal{F}$ with the subset $\{\diag(I,A)\in G: A\in \mathcal{F}\}$ of $G$.

Thus, as a corollary to Proposition \ref{isom}, every $g\in G$ is in the same orbit as $p = \diag(I,A)\in G$ where $I$ is the $(n-2)\times (n-2)$ identity matrix, $A\in \mathcal{F}$ and where the last column of $A$ is $(|v_0(g)|, |v_n(g)|, |v_{n+1}(g)|)^t$.

In other words, the orbit of the set $\mathcal{F}\subseteq G\times G$ under the $\Delta G\times K\times N$ action is all of $G\times G$.  It follows easily that the orbit of a dense subset of $\mathcal{F}$ is dense in $G\times G$.  Since $\pi:G\times G\rightarrow \Delta G\backslash G\times G/U$ is a submersion, $\pi$ maps open dense sets to open dense sets.  It follows that if we show the set of points in $\mathcal{F}$ which project to positively curved points in $G/U$ is dense, that $G/U$ is almost positively curved.  We summarize this in the following proposition.

\begin{proposition}\label{opden}
Consider the set of points $p\in \mathcal{F}$ which project to points in $G/U$ for which every $2$-plane has positive sectional curvature.  If this set is dense in $\mathcal{F}$, then $G/U$ is almost positively curved.

\end{proposition}

To actually compute, we use the following paramaterization of points in $\mathcal{F}$.

\begin{lemma}\label{lem:so3form}  Suppose $A \in \mathcal{F}$.  Then there are unique $\theta,\alpha\in [0,\pi/2]$ with $$A = \begin{bmatrix} \cos\theta & 0 & \sin \theta\\ -\cos\alpha \sin \theta & \sin \alpha & \cos\alpha \cos\theta\\ -\sin\alpha \sin\theta  & -\cos\alpha & \sin\alpha \cos\theta\end{bmatrix}.$$

\end{lemma}

\begin{proof}
Because the first row has length one and the last entry of the first row is non-negative, the first row has the form $(\pm \cos\theta, 0, \sin\theta)$ for a unique $\theta\in [0,\pi/2]$.  Orthogonality of the last two columns, together with the fact that each entry in the last column of $A$ and $A_{22}$ are all non-negative, implies $A_{32}\leq 0$.  Thus, the middle column of $A$ has the form $(0,\sin\alpha, -\cos\alpha)^t$ for some unique $\alpha\in [0,\pi/2]$.

 Now the form of last column of $A$ is determined using the fact that the entries are non-negative, and that is has unit length and is orthogonal to the second column.  Specifically, since it is unit length, we have $A_{23}^2 + A_{33}^2 = \cos^2\theta$, so $A_{23} = \cos\theta \cos \eta$ and $A_{33} = \cos\theta \sin\eta$ for some $\eta \in [0,2\pi)$.  Non-negativity then forces $\eta \leq \pi/2$.  Orthogonality with the second column shows $\tan \alpha = \tan \eta$; the bounds on $\alpha$ and $\eta$ now imply $\eta = \alpha$, as claimed.

Finally, the cross product of the second and third columns gives the first.  In particular, $A_{11} = \cos\theta$.

\end{proof}

Finally, in order to use Proposition \ref{eqns}, we must argue that both $G/K$ and $K/U$ are positively curved if $G$ is given a bi-invariant metric.  For $G/K \cong \mathbb{H}P^n$, the only $G$-invariant metric is, up to scaling, the Fubini-Study metric, so is positively curved.  On the other hand, $K/U \cong S^{4n-1}$ admits many $K$-invariant metrics, and the normal homogeneous metric is not the round metric on $S^{4n-1}$.  Nonetheless, the bi-invariant metric on $G$ restricts to a bi-invariant metric on $K$, and the following lemma shows the induced metric on $K/U\cong S^{4n-1}$ is positively curved.  Hence Proposition \ref{eqns} applies to all of these spaces.

\begin{lemma}

The bi-invariant metric $\langle \cdot ,\cdot \rangle_0$ on $Sp(n)$ induces a positively curved metric on $Sp(n)/Sp(n-1) \cong S^{4n-1}$.

\end{lemma}

\begin{proof}  It is well known (see, e.g., \cite[Corollary 3.33]{CE}) that the curvature of a normal homogeneous space $G/U$ is given by $\sec(X,Y) = \frac{1}{4}\| [X,Y]_\mathfrak{q}\|^2 + \| [X,Y]_\mathfrak{u}\|^2$ with $\mathfrak{g} = \mathfrak{u}\oplus \mathfrak{q}$.  In particular, $\sec(X,Y) = 0$ iff $[X,Y]=0$.  

Now, suppose for a contradiction that $\sigma \subseteq T_{eSp(n-1)} Sp(n)/Sp(n-1)$ is a $2$-plane with zero sectional curvature, where $Sp(n-1)$ is embedded into $Sp(n)$ as top left $(n-1)\times(n-1)$ block.  We let $\{X,Y\}$ denote a basis of $\sigma$.  We may interpret $X,Y\in \mathfrak{q} = \mathfrak{sp}(n-1)^\bot\subseteq \mathfrak{sp}(n)$.  Since the adjoint action of $Sp(n-1)$ on $\mathfrak{q}$ splits as a sum of the standard representation (which acts transitively on the unit sphere) and three trivial representations, we may assume without loss of generality that $X$ has the form $X = \begin{bmatrix}  0 & 0& ... & 0 & x_1\\ 0 & 0 & ... &0 & 0\\ \vdots & \vdots & \ddots & \vdots & \vdots\\ 0 & 0 & ... &0 & 0\\ -x_1 & 0 & ... & 0 & x_n\end{bmatrix}$ with $x_1 \in \mathbb{R}$ and $x_n \in \operatorname{Im}(\mathbb{H})$.

Now, the action by any matrix in $\{1\}\times Sp(n-2)\subseteq Sp(n-1)$ fixes $X$.  Using this action, we may assume without loss of generality that $Y$ has the form $Y = \begin{bmatrix}0 & 0& 0 & ... & 0 & y_1\\ 0 & 0 & 0 &  ... &0 & y_2\\ 0 & 0 & 0 & ... & 0 & 0\\ \vdots & \vdots & \vdots & \ddots & \vdots & \vdots \\ 0 & 0 & 0 & ... &0 & 0\\ -\overline{y}_1 & -y_2 & 0& ... & 0 & y_n\end{bmatrix}$ with $y_1\in\mathbb{H}$, $y_2 \in \mathbb{R}$ and $y_n \in \operatorname{Im}(\mathbb{H})$.

Now we compute $[X,Y] = XY-YX$ to be $$\begin{bmatrix} -x_1 \overline{y}_1  & -x_1 y_2 & 0 & ... & x_1 y_n\\ 0 & 0 & 0 & ... & 0\\  \vdots & \vdots & \vdots & \ddots & \vdots\\ -x_n \overline{y}_1 & -x_n y_2 & 0 & ... & -x_1 y_1 + x_n y_n\end{bmatrix} - \begin{bmatrix} -y_1 x_1 & 0 & ... & y_1 x_n\\ -y_2 x_1 & 0 & ... & y_2 x_n \\ \vdots & \vdots & \ddots & \vdots\\ -x_1 y_n & 0 & ... & -\overline{y}_1 x_1 + y_n x_n  \end{bmatrix}  $$ so vanishes iff \begin{align*}  2x_1\operatorname{Im}(y_1) &= 0 \\ x_1 y_2 &= 0\\ x_1 y_n -  y_1x_n &= 0\\ y_2 x_n &=0 \\ -2x_1\operatorname{Im}(y_1) + [x_n,y_n] &=0.\end{align*}

Assume initially that $x_n = 0$, so $x_1\neq 0$.  Then the first equation implies $y_1\in \mathbb{R}$.  The second equation implies $y_2 = 0$ and the third implies $y_n = 0$.  Thus $X$ and $Y$ are linearly dependent, giving a contradiction.

Thus, we must have $x_n\neq 0$.  The first and fifth equations taken together imply that $x_n$ and $y_n$ are linearly dependent over $\mathbb{R}$.  By subtracting an appropriate multiple of $X$ from $Y$, we may assume $y_n = 0$.  Then the third equation implies $y_1 = 0$.  Finally, the fourth equation gives $y_2 = 0$ so $Y = 0$ and $\{X,Y\}$ is not linearly independent, a contradiction.

\end{proof}

\section{\texorpdfstring{Almost positive curvature on $M_2$ and the two circle quotients.}{Almost positive curvature on M2 and the two circle quotients}}\label{lowdimalmpos}

In this section, we show the metrics constructed in Section \ref{metrics} are almost positively curved in the case of $M_2 = Sp(3)/Sp(1)\times Sp(1)$ and the two circle quotients $R_2 = \Delta S^1\backslash Sp(3)/Sp(1)^2$ and $Q_2 = Sp(3)/Sp(1)^2 S^1$.  From O'Neill's formulas, it is enough to show that $M_2$ is almost positively curved.

We denote $G = Sp(3)$, $K = Sp(2)\times Sp(1)$, and $U = Sp(1)\times \{1\}\times Sp(1)\subseteq K$, with Lie algebras $\mathfrak{g} = \mathfrak{sp}(3)$, etc.  Let $p\in \mathcal{F}$.  By Proposition \ref{opden}, $M_2$ is almost positively curved if the set of points in $\mathcal{F}$ for which all 2-planes are positively curved is dense in $\mathcal{F}$.

So, assume $[(p,e)]\in G/U \cong \Delta G\backslash G\times G/U$ has at least one zero-curvature plane.  As verified in Section \ref{metrics}, the metrics on $G/K \cong \mathbb{H}P^2$ and $K/U \cong S^7$ satisfy the hypothesis of Proposition \ref{eqns}, so we assume $X = X_\mathfrak{k}$ and $Y = Y_\mathfrak{p}$ are linearly independent vectors in $\mathfrak{g}$ which satisfy all the conditions of Proposition \ref{eqns}.

A simple calculation shows that $X = X_\mathfrak{k}$ satisfying \ref{con1} has the form $$X = \begin{bmatrix} 0 & a & 0\\ -\overline{a} & b & 0\\ 0 & 0 &  0 \end{bmatrix}$$ where $a\in \mathbb{H}$ and $b \in \operatorname{Im}\mathbb{H}$.  Likewise, since $Y\in \mathfrak{p}$, $Y$ has the form $$Y = \begin{bmatrix} 0 & 0& c \\ 0 & 0 & d\\ -\overline{c} & -\overline{d} & 0\end{bmatrix}$$ with $c,d\in \mathbb{H}$.

The form of $X$ and $Y$ are further constrained by \ref{con2}.

\begin{proposition}\label{niceform}  The vectors $X$ and $Y$ satisfy \ref{con2} iff $a=d=0$.
\end{proposition}

\begin{proof}  We compute $0=[X,Y]= \begin{bmatrix}0 &0 & ad \\ 0 & 0 & -\overline{a}c +bd \\ -\overline{d}\overline{a} & \overline{c}a + \overline{d}b &0  \end{bmatrix}$ which vanishes iff $\begin{bmatrix} ad\\ -\overline{a}c + bd \end{bmatrix} = 0.$  If $a\neq 0$, then the first entry forces $d = 0$, and then the second entry forces $c = 0$, that is, $Y = 0$.  Since $\{X,Y\}$ is linearly independent, this is a contradiction, so we must have $a = 0$, and thus, $b\neq 0$.  Then the second entry gives $bd = 0$, so $d=0$.  
\end{proof}

To apply condition \ref{con3}, we first compute $Ad_{p^{-1}} X = p^{-1} X p$ and $Ad_{p^{-1}}Y = p^{-1}Yp$.  We recall we are assuming $p\in \mathcal{F}$, so $p$ is a $3\times 3$ matrix with entries as in Lemma \ref{lem:so3form}.  To aid the calculation, we note the entries of $p$ are real, so commute with the entries of $X$ and $Y$.  Then a simple calculation shows that $$Ad_{p^{-1}} X = b\begin{bmatrix} \cos^2\alpha \sin^2\theta & -\cos\alpha \sin\alpha \sin\theta & -\cos^2\alpha \sin\theta \cos\theta\\ -\cos\alpha \sin\alpha \sin\theta & \sin^2\alpha & \cos\alpha\sin\alpha \cos\theta\\ -\cos^2\alpha \cos\theta\sin\theta & \cos \alpha \sin\alpha \cos\theta & \cos^2\alpha \cos^2\theta\end{bmatrix}$$ and that $Ad_{p^{-1}} Y$ is given by $$\begin{bmatrix} \overline{c} \sin\alpha \cos\theta\sin\theta - c\sin\alpha\cos\theta\sin\theta & -c\cos\alpha \cos\theta & \overline{c} \sin\alpha \sin^2\theta + c \sin\alpha \cos^2\theta \\ \overline{c} \cos\alpha \cos\theta & 0 & \overline{c} \cos\alpha\sin\theta\\ -c \sin\alpha \sin^2\theta -\overline{c} \sin\alpha\cos^2\theta & -c\cos\alpha\sin\theta & -\overline{c} \sin\alpha \cos\theta\sin\theta + c\sin\alpha \cos\theta\sin\theta\end{bmatrix}.$$

So, \ref{con3} is satisfied iff $$V := \begin{bmatrix} -b\cos^2 \alpha \cos\theta \sin\theta \\ b\cos\alpha \sin\alpha \cos\theta\end{bmatrix} \text{ and } W :=\begin{bmatrix} \overline{c} \sin\alpha \sin^2\theta + c \sin\alpha \cos^2\theta \\ \overline{c} \cos\alpha \sin\theta \end{bmatrix}$$ are linearly dependent over $\mathbb{R}$.

Recalling that $\alpha,\theta\in [0,\pi/2]$, we note that $V$ is identically zero for some non-zero $b\in \operatorname{Im}\mathbb{H}$ iff $\alpha = \pi/$2, or $\theta = \pi/2$, or $\theta = \alpha = 0$.  Clearly, if one of the conditions is satisfied, then, from Proposition \ref{eqns}, there are zero-curvature planes at $p$.  Similarly, $W$ is identically zero for some non-zero $c\in \mathbb{H}$ iff $\alpha = \pi/2$ and $\theta = \pi/4$, or if $\alpha = \theta = 0$, and again, there will be zero-curvature planes at a point $p$ satisfying one of these conditions.  For the remainder of this section, we assume that $\theta,\alpha\in(0,\pi/2)$ and that $\theta \neq \pi/4$, so that, in particular, $V$ and $W$ are non-zero vectors for any non-zero choices of $b$ and $c$.  Clearly, there is an open dense subset $\mathcal{F}_1\subseteq \mathcal{F}$ for which this condition on $\theta$ and $\alpha$ holds.  Because $V$ and $W$ are non-zero, \ref{con3} is satisfied iff $V = W$ for some $0\neq b\in \operatorname{Im}\mathbb{H}$ and $0\neq c\in \mathbb{H}$.

We are now in a position to show that $M_2$ has an open dense set of points for which all $2$-planes are positively curved.

\begin{proposition}\label{case1proof}

Suppose there are non-zero $b\in \operatorname{Im}\mathbb{H}$ and $c\in \mathbb{H}$ for which $V = W$.  Then $\tan^2\alpha = \frac{\sin^2\theta}{\cos^2\theta - \sin^2\theta}$

\end{proposition}

\begin{proof}

First, we rewrite the first entry of $W$ as $\sin\alpha( \operatorname{Re}(c) + (\cos^2\theta - \sin^2\theta) \operatorname{Im}(c))$.  Since $b\in \operatorname{Im}\mathbb{H}$, the equation $$-\cos^2\alpha \cos\theta \sin\theta b = \sin\alpha( \operatorname{Re}(c) + (\cos^2\theta - \sin^2\theta) \operatorname{Im}(c)),$$ which comes from the first component of the equation $V=W$, implies that $c$ is purely imaginary.  So we may now rewrite the first entry of $W$ as $\sin\alpha(\cos^2 \theta - \sin^2\theta) c$.  Thus, we see $$-\cos^2\alpha \cos\theta \sin\theta b = \sin\alpha(\cos^2\theta - \sin^2\theta)c,$$ so $$b = \frac{\sin\alpha(\cos^2\theta - \sin^2\theta)}{-\cos^2\alpha \cos\theta \sin\theta } .$$

Substituting this into the second component of the equation $V = W$, noting that $y = \operatorname{Im}(y)$ implies $\overline{y} = -y$, and canceling $c$, we obtain the equation $$ \cos\alpha\sin\alpha \cos\theta \frac{\sin\alpha(\cos^2\theta - \sin^2\theta)}{-\cos^2\alpha \cos\theta\sin\theta} = -\cos\alpha \sin\theta.$$  This simplifies to $\tan^2\alpha = \frac{\sin^2\theta}{\cos^2\theta - \sin^2\theta}$, as claimed. 
\end{proof}

If we let $\mathcal{F}_2\subseteq \mathcal{F}_1$ be the subset of of $\mathcal{F}_1$ with $\tan^2\alpha \neq \frac{\sin^2\theta}{\cos^2\theta - \sin^2\theta}$, then $\mathcal{F}_2$ consists of points $p\in \mathcal{F}$ which project to points in $M_2$ for which every $2$-plane has positive sectional curvature.  Clearly, $\mathcal{F}_2$ is an open dense subset of both $\mathcal{F}_1$ and $\mathcal{F}$.  Thus, from Proposition \ref{opden}, $M_2$ is almost positively curved.  This completes the proof of Theorem \ref{main}.

\section{\texorpdfstring{Open sets of zero-curvature points on $M_n$ for $n\geq 3$}{Open sets of zero-curvature points on Mn for n >=3}}\label{hidim}

In this section, we show Tapp's metrics \cite{Ta1} are not always almost positively curved.  More specifically, we show that Tapp's quasi-positively curved metrics on $M_n = Sp(n+1)/Sp(n-1)Sp(1)$ with $n\geq 3$ are not almost positively curved.   We let $G = Sp(n+1)$ and $U = Sp(n-1)\times Sp(1)$, with $U$ embedded in $G$ as $(A,q)\mapsto \diag(A,1,q)$.

In \cite{Ta1}, Tapp shows his metrics are, up to scaling, isometric to those defined in Section \ref{metrics}, with $K = Sp(n)\times Sp(1)$.  Proposition \ref{isom} applies in this case, so every point in $\Delta G\backslash G\times G/U$ is isometrically equivalent to a point in $\mathcal{F}$.

We now find an open subset of $G/U$ for which every point has infinitely many zero-curvature planes.

To that end, given $A \in \mathcal{F}$ (so $A$ has the form given by Lemma \ref{lem:so3form}), we make the following definitions:

$$ \mu = \sqrt{\tan^2\theta \csc^2\alpha - 1} \text{ and } \eta = \frac{1}{\mu} \frac{\sin\theta (\cos^2\theta - \sin^2\theta)}{\cos \alpha \sin^2\alpha \cos^3\theta}$$

We let $Z\subseteq \mathcal{F}$ denote the open set of points for which $\mu > 0$ and $\eta$ is defined, that is, where the denominator of $\eta$ is non-zero.  Since the orbit through a point $(g_1,g_2)\in G\times G$ under the natural $\Delta G\times K\times N$ action is determined by the lengths $|v_i(g_2^{-1} g_1)|$ for $i = 0,n,n+1$ (Proposition \ref{isom}), we see that the set of points in $G\times G$ whose orbits pass through $Z\times \{1\}$ is open in $G\times G$.  In particular, if we can show that for every $p = \diag(I,A)$ with $A\in Z$, that the point $[(p,e)]\in \Delta G\backslash G\times G/U$ has at least one zero-curvature plane, then Theorem \ref{main2} must be true.

\begin{proposition}\label{notalmost1}

Suppose $A\in Z$ and $p = \diag(I,A)\in G$.  Then there are infinitely many zero-curvature planes at the point $[(p,e)]\in \Delta G\backslash G\times G/U\cong  G/U$.

\end{proposition}

\begin{proof}  Fix any purely imaginary unit length quaternion $b$.  We set $X\in \mathfrak{g} = \mathfrak{sp}(n+1)$ to the matrix which is zero everywhere except the bottom right $4\times 4$ block, where it is $X_0:=\begin{bmatrix} 0 & 0 & 1 & 0 \\ 0 & 0 & \mu b & 0\\ -1 & \mu b & \eta b & 0 \\0 & 0 & 0 & 0\end{bmatrix}$.  Note that $X = X_{\mathfrak{k}}$.

Likewise, we define $Y\in \mathfrak{g}$ to be the matrix which is zero everywhere except the bottom right $4\times 4$ block, where it is $Y_0:=\begin{bmatrix}0 & 0 & 0 & 1\\ 0 & 0 & 0 & -\frac{b}{\mu}\\ 0 & 0 & 0 & 0 \\ -1 &  -\frac{b}{\mu} & 0 & 0 \end{bmatrix}$.  Note that $Y = Y_\mathfrak{p}$.  We claim that $X$ and $Y$ satisfy all the conclusions of Proposition \ref{eqns}, so there is a zero curvature plane at $[(p,e)]\in G/U$.

Clearly, both $X$ and $Y$ are orthogonal to $\mathfrak{u}$, so we may focus on \ref{con2}.  Because of the block form of $X$ and $Y$, we see that $[X,Y]=0$ iff $[X_0, Y_0] = 0$.  Computing the latter, the only potentially non-zero entries are $[X_0,Y_0]_{3,4} = -1-b^2$ and $[X_0,Y_0]_{4,3} = 1+b^2$.  Since $b$ is a purely imaginary unit length quaternion, $b^2 = -1$, so \ref{con2} is satisfied.

In order to verify \ref{con3}, we compute $Ad_{p^{-1}} X$ and $Ad_{p^{-1}}Y$.  Due to the block form of $X,Y,$ and $p$, it follows that outside of the bottom right $4\times 4$ block, every entry of both $Ad_{p^{-1}}X$ and $Ad_{p^{-1}}Y$ vanishes.  Further, the bottom right $4\times 4$ block of $Ad_{p^{-1}}X$ is equal to $Ad_{\diag(1,A)^{-1}} X_0$, and likewise for $Y$.  A simple calculation now gives the last column of $Ad_{\diag(1,A)^{-1}}X_0$ and $Ad_{\diag(1,A)^{-1}}Y_0$ as $$\begin{bmatrix}\cos \alpha cos\theta \\ \mu b \cos\alpha(\cos^2\theta - \sin^2\theta) - \eta b \cos^2\alpha \cos\theta\sin\theta \\ \mu b\sin\alpha \sin\theta + \eta b \cos\alpha \sin\alpha \cos\theta\\ \mu b \cos\alpha\cos\theta\sin\theta + \eta b \cos^2\alpha \cos^2\theta  \end{bmatrix} \text{ and } \begin{bmatrix}\sin\alpha \cos\theta \\ \frac{b}{\mu} \sin\alpha(\sin^2\theta - \cos^2\theta)\\ \frac{b}{\mu} \cos \alpha \sin\theta \\ -2\frac{b}{\mu}\sin \alpha \cos\theta\sin\theta  \end{bmatrix}$$ respectively.  In particular, the non-zero entries of the $\mathfrak{p}$ components of $Ad_{p^{-1}}X$ and $Ad_{p^{-1}} Y$ can be identified with the vectors $$V:=\begin{bmatrix}\cos \alpha cos\theta \\ \mu b \cos\alpha(\cos^2\theta - \sin^2\theta) - \eta b \cos^2\alpha \cos\theta\sin\theta \\ \mu b\sin\alpha \sin\theta + \eta b \cos\alpha \sin\alpha \cos\theta\end{bmatrix}$$ and $$W:=\begin{bmatrix} \sin\alpha \cos\theta \\ \frac{b}{\mu} \sin\alpha(\sin^2\theta - \cos^2\theta)\\ \frac{b}{\mu} \cos \alpha \sin\theta \end{bmatrix}$$ respectively.

Now, \ref{con3} is verified iff $V$ and $W$ are linearly dependent.  The fact that $\mu$ and $\eta$ are defined means that the first entry of both $V$ and $W$ is non-zero, so \ref{con3} is verified iff $V = \cot\alpha W$.  Since $\mu$ and $b$ are both non-zero, $V = \cot\alpha W$ iff $\frac{\mu}{b} V = \frac{\mu}{b}\cot\alpha W$.  The second and third entries of the equation $\frac{\mu}{b} V = \frac{\mu}{b}\cot\alpha W$ are $$\left\{ \begin{array}{lll} \mu^2  \cos\alpha(\cos^2\theta - \sin^2\theta) &- \eta\mu  \cos^2\alpha \cos\theta\sin\theta  &=  \cos\alpha(\sin^2\theta - \cos^2\theta)\\ \mu^2 \sin\alpha \sin\theta &+ \eta\mu  \cos\alpha \sin\alpha \cos\theta &=  \frac{\cos^2 \alpha \sin\theta}{\sin\alpha} \end{array} \right. .$$

According to Lemma \ref{computation} below, $V$ and $W$ are linearly dependent iff $$\mu^2 = \tan^2\theta \csc^2\alpha - 1 \text{ and } \eta\mu = \frac{\sin\theta (2\cos^2\theta - 1)}{\cos \alpha \sin^2\alpha \cos^3\theta}.$$  From the definition of $\mu$ and $\eta$ above, $V$ and $W$ are linearly dependent, so $[(p,e)]\in \Delta G\backslash G\times G/U$ has zero-curvature planes.

Since $b$ was an arbitrary unit length quaternion, this gives infinitely many zero-curvature planes at $p\in Z$. 

\end{proof}

So, establishing the following Lemma completes the proof of Theorem \ref{main2}.

\begin{lemma}\label{computation}The solution to the system $$\left\{ \begin{array}{lll} \mu^2  \cos\alpha(\cos^2\theta - \sin^2\theta) &- \eta\mu  \cos^2\alpha \cos\theta\sin\theta  &=  \cos\alpha(\sin^2\theta - \cos^2\theta)\\ \mu^2 \sin\alpha \sin\theta &+ \eta\mu  \cos\alpha \sin\alpha \cos\theta &=  \frac{\cos^2 \alpha \sin\theta}{\sin\alpha} \end{array} \right.$$ is given by $\mu^2 = \tan^2\theta \csc^2\alpha - 1$ and $\eta\mu = \frac{\sin\theta (2\cos^2\theta - 1)}{\cos \alpha \sin^2\alpha \cos^3\theta}.$

\end{lemma}

\begin{proof}  Viewing the system as a linear system in the variables $\mu^2$ and $\eta\mu$, we solve via Cramer's rule.  The denominator is given by \begin{align*} &\cos \alpha (\cos^2\theta - \sin^2\theta) \cos\alpha \sin \alpha \cos\theta + \sin\alpha \sin\theta  \cos^2\alpha \cos\theta \sin\theta\\ =& \cos^2\alpha \sin\alpha \cos\theta (\cos^2\theta - \sin^2\theta  + \sin^2\theta)\\ =& \cos^2\alpha \sin\alpha \cos^3\theta. \end{align*}

For $\mu^2$, the numerator is given by \begin{align*}&\cos \alpha (\sin^2\theta - \cos^2\theta) \cos\alpha \sin\alpha \cos\theta + \frac{\cos^2\alpha \sin\theta}{\sin \alpha} \cos^2\alpha \cos\theta \sin \theta \\ =& \cos^2\alpha \sin\alpha \cos\theta(\sin^2\theta(1+\cot^2 \alpha) - \cos^2\theta)\\ =& \cos^2\alpha \sin\alpha \cos\theta (\sin^2\theta \csc^2 \alpha - \cos^2\theta).\end{align*}

Thus, \begin{align*} \mu^2 &= \frac{\cos^2\alpha \sin\alpha \cos\theta (\sin^2\theta \csc^2 \alpha - \cos^2\theta)}{\cos^2\alpha \sin\alpha \cos^3\theta}\\ &= \tan^2\theta \csc^2\alpha - 1\end{align*}

Similarly, the numerator of $\eta\mu$ is given by \begin{align*}&\cos\alpha(\cos^2\theta - \sin^2\theta) \frac{\cos^2\alpha \sin\theta}{\sin \alpha} - \sin\alpha \sin\theta \cos\alpha(\sin^2\theta - \cos^2\theta)\\ =& \cos\alpha \sin\theta (\cos^2\theta - \sin^2\theta)\left(\frac{\cos^2\alpha}{\sin\alpha} + \frac{\sin^2 \alpha}{\sin\alpha} \right)\\ =& \cot\alpha \sin\theta(\cos^2\theta - \sin^2\theta).\end{align*}

Thus, \begin{align*} \eta\mu &= \frac{\cot\alpha \sin\theta(\cos^2\theta - \sin^2\theta)}{\cos^2\alpha \sin\alpha \cos^3\theta}\\ &= \frac{\sin\theta (\cos^2\theta - \sin^2\theta)}{\cos\alpha \sin^2\alpha\cos^3\theta}\end{align*}

\end{proof}

\section{\texorpdfstring{The topology of $M_2$, $Q_2$, and $R_2$}{The topology of M2, Q2, and R2}}\label{top}

We now compute the cohomology rings and characteristic class of $M_2 = Sp(3)/Sp(1)^2$ and the two circle quotients $R_2 = Sp(3)/Sp(1)^2 S^1$ and $Q_2 = \Delta S^1\backslash Sp(3)/Sp(1)^2$.  Singhof and Wemmer \cite{SiWe1} have shown $Sp(3)/Sp(1)^2$ is parallelizable.  To compute the cohomology ring of the homogeneous space $Sp(3)/Sp(1)^2$, we consider the chain of subgroups $Sp(1)^2\rightarrow Sp(2)\rightarrow Sp(3)$, where the embedding $Sp(2)\rightarrow Sp(3)$ is given by $\begin{bmatrix} a&b\\ c&d\end{bmatrix} \mapsto \begin{bmatrix} a& 0 & b\\ 0 & 1 & 0 \\ c& 0 & d\end{bmatrix}$.  The one has the homogeneous fibration $$S^4 \cong Sp(2)/Sp(1)^2\rightarrow Sp(3)/Sp(1)^2 \rightarrow Sp(3)/Sp(2) \cong S^{11}$$ showing $Sp(3)/Sp(1)^2$ is an $S^4$ bundle over $S^{11}$.  The Gysin sequence associated to this fiber bundle and Poincar\'e duality then imply the integral cohomology ring is isomorphic to that of $S^4\times S^{11}$.

As Kamerich \cite{Kam} showed in his thesis, $Sp(3)/Sp(1)^2$ and $S^4\times S^{11}$ are not homotopy equivalent.  We provide a short proof below for the convenience of the reader.

\begin{proposition}\label{nothe}  The homotopy groups $\pi_{10}(S^4\times S^{11})$ and $\pi_{10}(M_2)$ are not isomorphic, so $S^4\times S^{11}$ and $M_2$ are not homotopy equivalent.
\end{proposition}

\begin{proof}  We first note that $\pi_{10}(Sp(3))$ is in the stable range, so is given by Bott periodicity.  Thus, $\pi_{10}(Sp(3)) = 0$.  Recall (\cite[pg. 399]{Ha1}) $\pi_{10}(S^4\times S^{11})\cong \pi_{10}(S^4)\cong \mathbb{Z}_{24}\oplus \mathbb{Z}_3$ and $\pi_9(S^3) = \mathbb{Z}_3$.

Now, a portion of the long exact sequence of homotopy groups associated to the fibration $Sp(1)^2\rightarrow Sp(3)\rightarrow Sp(3)/ Sp(1)^2$ is $$ 0 = \pi_{10}(Sp(3))\rightarrow \pi_{10}(Sp(3)/Sp(1)^2)\rightarrow \pi_9(Sp(1)^2)\rightarrow ...$$  But $\mathbb{Z}_{24}\oplus\mathbb{Z}_3$ cannot inject into $\mathbb{Z}_3\oplus\mathbb{Z}_3$. 

\end{proof}

We now show $Q_2$ and $R_2$ have isomorphic cohomology rings, but that their first Pontryagin class mod $24$ are different.  Since this is a homotopy invariant \cite{AH}, this implies $Q_2$ and $R_2$ are homotopically distinct.

To do this, we first view both as biquotients in the form $H_i \backslash G/ K_i$:  $$Q_2 = \{e\}\backslash Sp(3)/Sp(1)^2\times S^1\text{ and } R_2 = \Delta S^1 \backslash Sp(3)/ Sp(1)^2$$ defined by two inclusions $H_i\times K_i\rightarrow G\times G$:

$$\text{for }Q_2,\text{ }(q_1,q_2,z)\rightarrow \left( I, \diag(q_1,z,q_2)\right)$$ and $$\text{for }R_2,\text{ } (q_1,q_2,z)\rightarrow \left(\diag(z,z,z), \diag(q_1,1,q_2)\right),$$ where $q_i\in Sp(1)$ and $z\in S^1$.

Letting $BG$ denote the classifying space of $G$, the quotient of a contractible space $EG$ by a free action of $G$, the inclusion $H_i\times K_i\rightarrow G\times G$ induces a map $BH_i\times BK_i\rightarrow BG\times BG.$  Using this map, Singhof \cite{Si1} proves the following theorem.

\begin{theorem}[Singhof]\label{fullrank} If the rank of $H\times K$ is equal to the rank of $G$, and if $H^\ast(BH), H^\ast(BK)$, and $H^\ast(BG)$ are all torsion free, then as algebras, $$H^\ast(G\bq (H\times K)) \cong H^\ast(BH)\otimes_{H^\ast(BG)} H^\ast(BK).$$
\end{theorem}

In order to determine the maps $H^\ast(BG)\rightarrow H^\ast(BH)$ and $H^\ast(BG)\rightarrow H^\ast(BK)$, we use a theorem of Borel \cite{Bo1}.

\begin{theorem}[Borel]\label{Borel} The inclusion map of a maximal torus $T\rightarrow Sp(n)$ induces an injective map $$H^\ast(BSp(n))\rightarrow H^\ast(BT)\cong \mathbb{Z}[x_1,...,x_n], \text{ with }|x_i| = 2$$ with image generated by the elementary symmetric polynomials in the squares of the $x_i$ variables.

\end{theorem}

Then, using the commutative diagram

\begin{diagram}
H^\ast(BG)& \rTo & H^\ast(BH)\\ 
\dTo & & \dTo\\ 
H^\ast(BT_G) & \rTo &H^\ast(BT_H)\\
\end{diagram}

induced from the natural inclusions, we compute the top map by computing the bottom map and restricting, and similarly for $H^\ast(BG)\rightarrow H^\ast(BK)$.

We now carry this out for the more difficult case of $R_2$.

We identify $H^\ast(BT_G)$ with $\mathbb{Z}[x_1, x_2, x_3]$ where $x_i\in H^2(BT_G)$ are the transgressions of the generators of the usual basis of $H^1(T_G)$ in the spectral sequence associated to $T_G\rightarrow ET_G\rightarrow BT_G$, and we similarly identify $H^\ast(BT_K) \cong \mathbb{Z}[y_1,y_2]$ and $H^\ast(BT_H)\cong \mathbb{Z}[u]$.

Then, Proposition \ref{Borel} identifies $H^\ast(BG)$ with $\mathbb{Z}[\sigma_1(x_i^2), \sigma_2(x_i^2), \sigma_3(x_i^2)]$ and identifies $H^\ast(BK)$ with $\mathbb{Z}[y_1^2, y_2^2]$.  Of course, since $H = T_H$, $H^\ast(BH)\cong H^\ast(BT_H) \cong \Z[u]$.

The map $H^\ast(BT_G)\rightarrow H^\ast(BT_K)$ maps $x_1$ to $y_1$, $x_2$ to $0$, and $x_3$ to $y_2$.  Thus, the map $H^\ast(BG)\rightarrow H^\ast(BK)$ is given as follows: $$\begin{matrix} \sigma_1(x_i^2) &=& x_1^2 + x_2^2 + x_3^2 &\mapsto& y_1^2 + y_2^2 \\ \sigma_2(x_i)^2 & = & x_1^2 x_2^2 + x_1^2 x_3^2 + x_2^2 x_3^2 &\mapsto& y_1^2 y_2^2\\ \sigma_3(x_i)^2 &=& x_1^2 x_2^2 x_3^2 & \mapsto & 0.\end{matrix}$$

 Similarly, the function $H^\ast(BT_G)\rightarrow H^\ast(BT_H)$ maps $x_i$ to $u$ for all $i$, and thus, the map $H^\ast(BG)\rightarrow H^\ast(BH)$ is given by $\sigma_1(x_i^2)\mapsto 3u^2$, $\sigma_2(x_i^2)\mapsto 3u^4$, and $\sigma_3(x_i^2) \mapsto u^6$.

Thus, Theorem \ref{fullrank} implies $$H^\ast(R_2)\cong \mathbb{Z}[y_1^2, y_2^2, u]/I $$ where $y_i$ and $u$ both have degree $2$ and $I$ the ideal generated by $y_1^2 + y_2^2  -3u^2$, $y_1^2 y_2^2 - 3u^4$, and $u^6$.  One sees easily that this is isomorphic to $\mathbb{Z}[y_1^2, u]/I_2$ where $I_2$ is generated by $3u^4 - 3 y_1^2 u^2 + y_1^4$ and $u^6$.

In a similar fashion, Theorem \ref{fullrank} can be used to show $$H^\ast(Q_2)\cong \mathbb{Z}[y_1^2, y_2^2, u]/J$$ where $J$ is the ideal generated by $y_1^2 + y_2^2 + u^2$, $y_1^2 y_2^2 + (y_1^2 + y_2^2)u^2 $, and $y_1^2 y_2^2 u^2 $, which is clearly isomorphic to $\mathbb{Z}[y_1^2, u]/J_2$ where $J_2$ is generated by $y_1^4 + y_1^2 u^2 + u^4$ and $u^6$.

\begin{proposition}\label{cohomring}The cohomology rings $\Z[y_1^2,u]/I_2$ and $\Z[y_1^2, u]/J_2$ are isomorphic.

\end{proposition}

\begin{proof}Consider the function $\phi:\Z[y_1^2, u]\rightarrow \Z[y_1^2,u]$ given by $\phi(u) = u$ and $\phi(y_1^2) = u^2 - y_1^2$.  It is easy to verify that $\phi^2$ is the identity function.

Also, since $\phi(u^6) = u^6 \in J_2$ and \begin{align*} \phi(3u^4 - 3 y_1^2 u^2 + y_1^4) &= 3u^4 - 3(u^2 - y_1^2)u^2 + (u^2 - y_1^2)^2\\ &= 3u^4 -3u^4 + 3y_1^2 u^2 + u^4 - 2u^2 y_1^2 + y_1^4\\ &= y_1^4 +  y_1^2u^2 + u^4\\ &\in J_2,\end{align*} $\phi(I_2)\subseteq J_2$, so $\phi$ induces a map from $\Z[y_1^2,u]/I_2$ to $\Z[y_1^2,u]/J_2$.  In a similar manner, it is easy to verify that $\phi(y_1^4 + y_1^2 u^2 + u^4) = y_1^4 - 3y_1^2 u^2 + 3u^4\in I_2$, so $\phi$ induces a map from $\Z[y_1^2, u]/J_2$ to $\Z[y_1^2,u]/I_2$.  Since $\phi^2$ is the identity, these induced maps are inverses of each other, so they are both isomorphisms.  

\end{proof}

We now set up notation in order to compute the first Pontryagin classes of $Q_2$ and $R_2$.  For $T\subseteq G$ a torus, we may use transgressions of generators of $H^1(T)$ in the spectral sequence $T\rightarrow ET\rightarrow BT$ as generators of $H^2(BT)$.  Since a weight of a representation of $G$ is an element of the weight lattice $\ker \exp$ with $\exp:\mathfrak{t}\rightarrow T$ the group exponential map, this allows us to interpret weights of a representation as elements of $\operatorname{Hom}(\ker \exp, \mathbb{Z}) \cong \operatorname{Hom}(\pi_1(T), \mathbb{Z})\cong H^1(T)$, which may then be interpreted, via transgressions, as elements of $H^2(BT)$.

Using this notation, Singhof \cite{Si1} proves the following theorem, adapted to the full rank case.

\begin{theorem}[Singhof]\label{char} Let $\Delta^+_G$ denote the set of positive roots of $G$, interpreted as elements of $H^2(BT_G)$ and similarly for $\Delta^+_H$ and $\Delta^+_K$.  Then the isomorphism in Theorem \ref{fullrank} gives an identification $$ p(H\backslash G/K) = \prod_{\beta \in \Delta^+_G}(1+\beta^2) \prod_{\gamma\in \Delta^+_H}(1+\gamma^2)^{-1} \prod_{\delta \in \Delta^+_K}(1+\delta^2)^{-1}.$$

\end{theorem}

We may now distinguish $Q_2$ and $R_2$.

\begin{proposition}\label{homotopytype}  The manifolds $Q_2$ and $R_2$ have distinct homotopy types.

\end{proposition}

\begin{proof}

In the notation of the previous computation of $H^\ast(R_2)$, the positive roots of $Sp(3)$ are $2x_i$ and $x_i \pm x_j$ for $1\leq i < j\leq 3$, while for $Sp(1)\times Sp(1)$, they are $2y_1$ and $2y_2$.  Of course, a circle $S^1$ has no positive roots.  Using the notation $(x_1\pm x_2)^2$ to mean $(x_1 + x_2)^2 + (x_1 - x_2)^2 = 2(x_1^2 + x_2^2)$, it follows from Theorem \ref{char} that \begin{align*} p_1 = & \sum \beta^2 - \sum \gamma^2 - \sum \delta^2 \\ =& 4(x_1^2 + x_2^2 + x_3^2) + (x_1\pm x_2)^2 + (x_1\pm x_3)^2+ (x_2\pm x_3)^2\\ &-4y_1^2 - 4y_2^2\\ =& 8(x_1^2 + x_2^2 + x_3^2) - 4(y_1^2 + y_2^2).\end{align*}

Now, via the inclusion $BH\rightarrow BG$, we see $\sigma_1(x_i^2) = x_1^2 + x_2^2 + x_3^2$ maps to $0$ for $Q_2$ and to $3z^2$ for $R_2$.  Since $y_1^2 + y_2^2 = -z^2$ in $H^\ast(Q_2)$ and $y_1^2 + y_2^2 = 3z^2$ in $H^\ast(R_2)$, we see that $p_1(Q_2) = 4z^2$ while $p_1(R_2) = 24z^2  -12z^2 = 12z^2$.  Now, one easily checks that $H^4(Q_2;\mathbb{Z}_{24})/p_1 = \mathbb{Z}_{24} \oplus \mathbb{Z}_6$ while $H^4(R_2;\mathbb{Z}_{24})/p_1 = \mathbb{Z}_{24}\oplus \mathbb{Z}_2$.  Thus, there is no isomorphism $H^4(Q_2; \mathbb{Z}_{24})\rightarrow H^4(R_2; \mathbb{Z}_{24})$ which preserves $p_1$.  Since this is a homotopy invariant \cite{AH}, it follows that $Q_2$ and $R_2$ are not homotopy equivalent. 

\end{proof}

\

The only previously known examples of simply connected almost positively curved manifolds in dimension $14$ and $15$ are due to Wilking \cite{Wi}.  In dimension $15$, they are $T^1 S^8$ and the homogeneous space $U(5)/U(3)S^1_{kl}$, while in dimension $14$, they are $\Delta SO(2)\backslash SO(9)/SO(7)$ and $P_\mathbb{C} T^1\mathbb{C}P^4$, the projectivized unit tangent bundle to $\mathbb{C}P^4$.

Now, $T^1 S^8$ is $6$-connected, while $\pi_2(U(5)/U(3)S^1_{kl}) \cong \Z$.  On the other hand, $\pi_4(M_2) \cong \Z$ while $\pi_2(M_2) = 0$, so the $15$ dimensional example is distinct up to homotopy from the previously known examples.

Further, $\Delta SO(2) \backslash SO(9)/SO(7)$ is a circle quotient of $SO(9)/SO(7)\cong T^1 S^8$, so again has $\pi_4$ trivial, while both $Q_2$ and $R_2$ have $\pi_4$ isomorphic to $\mathbb{Z}$.  Finally, $P_\mathbb{C} T^1 \mathbb{C}P^4$ fits into a fiber bundle $S^1\rightarrow T^1\mathbb{C}P^4\rightarrow P_\mathbb{C}T^1\mathbb{C}P^4$, so $\pi_2(P_\mathbb{C}T^1 \mathbb{C}P^4)\cong \Z^2$, while $\pi_2(Q_2) \cong \pi_2(R_2)\cong \Z$.  So the two $14$-dimensional examples are distinct up to homotopy from the previously known examples as well.


\bibliographystyle{plain}

\end{document}